\documentclass[12pt]{article}
\usepackage[letterpaper,text={6.5in,9in}]{geometry}

\usepackage{amsfonts}
\usepackage{amssymb}
\usepackage{amsmath}

\usepackage{xparse}  

\usepackage{color}
\usepackage[colorlinks=true,citecolor=black,linkcolor=black,urlcolor=blue]{hyperref}
\usepackage{graphicx}

\usepackage{enumerate}

\allowdisplaybreaks

\newcommand{\ttt}{Tic-Tac-Toe}
\newcommand{\MB}{Maker--Breaker}
\newcommand{\M}{Maker}
\newcommand{\B}{Breaker}
\newcommand{\cH}{{\cal H}}

\newcommand{\naturals}{{\mathbb N}}

\newcommand{\cQ}{{\cal Q}}

\newcommand{\PS}{PS(4,2)}
\newcommand{\ps}{p}

\newcommand{\binSetsPS}[1][n]{\mbox{\em BinProds}(4 \times #1)}

\NewDocumentCommand{\binPS}{O{n} O{4}}{\mbox{\em Bin}PS(#2 \times #1)}

\NewDocumentCommand{\binPSR}{O{4} O{n}}{\mbox{\em Bin}PS^R(#1 \times #2)}

\newcommand{\Rs}{R(s)}

\NewDocumentCommand{\binPSP}{O{n} O{4}  m  m }{\mbox{\em Bin}PS_{#3}^{R (#4)}(#2 \times #1)}

\newcommand{\binSetsPSP}[2]{\mbox{\em BinSets}PS_{#1}^{R(#2)}(4,2)}

\newcommand{\smallBinPSP}[2]{\mbox{\em Bin}PS_{#1}^{R(#2)}(3,2)}
\newcommand{\smallBinSetsPSP}[2]{\mbox{\em BinSets}PS_{#1}^{R(#2)}(3,2)}

\newcommand{\BPS}{Breaker's win pairing strategy}
\newcommand{\BPSs}{Breaker's win pairing strategies}

\newcommand{\binform}{product-set}   
\newcommand{\binforms}{product-sets}
\newcommand{\binset}{bin-set}
\newcommand{\binsets}{bin-sets}

\newcommand{\block}{block}
\newcommand{\blocks}{blocks}
\newcommand{\blocked}{blocked}

\newcommand{\LD}{LD}
\newcommand{\ld}{lower dimensional}
\newcommand{\HD}{HD}
\newcommand{\hd}{higher dimensional}

\newcommand{\Mbin}{{\bf M}-bin}

\newcommand{\Star}{star}   
\newcommand{\Stars}{stars}	

\newcommand{\bodd}{{\bf 1}}
\newcommand{\beven}{{\bf 0}}
\newcommand{\bv}{{\bf P}}      
\newcommand{\bA}{{\bf a}}

\newcommand{\bva}{{\bf P}^{(0)}}   
\newcommand{\bvb}{{\bf P}^{(1)}}     
\newcommand{\bvc}{{\bf P}^{(2)}}	
\newcommand{\bvd}{{\bf P}^{(3)}}	
\newcommand{\bvj}{{\bf P}^{(j)}}	
\newcommand{\bvP}[1]{{\bf P}^{(#1)}}	

\newcommand{\rS}[1]{S|_{\mbox{\scriptsize bin }#1}}
\newcommand{\re}[1]{\vec{e}|_{\mbox{\scriptsize bin }#1}}
\newcommand{\rv}[1]{\vec{v}|_{\mbox{\scriptsize bin }#1}}
\newcommand{\rSet}[2]{#1|_{\mbox{\scriptsize bin }#2}}

\newcommand{\bevenP}[1]{{\bf 0}^{(#1)}}
\newcommand{\boddP}[1]{{\bf 1}^{(#1)}}

\newcommand{\Index}{\mbox{Index}}

\newcommand{\sundberg}{Sundberg}

\newcommand{\ex}{\mbox{ex}}

\definecolor{darkgreen}{rgb}{0, .5, 0}
\definecolor{darkmagenta}{rgb}{.55,0,.55}

\newtheorem{theorem}{Theorem}
\newtheorem{lemma}{Lemma}

\newtheorem{corollary}{Corollary}

\begin{document}

\title{Pairing strategies for the \MB\ game on the hypercube with subcubes as
winning sets
}

\author{Ramin Naimi \qquad Eric Sundberg \\
\small Mathematics Department \\[-0.8ex]
\small Occidental College \\[-0.8ex]
\small Los Angeles, CA \\[-0.8ex]
\small\tt \{rnaimi, sundberg\}@oxy.edu \\
} 

\maketitle


\begin{abstract}
We consider the \MB\ positional game on the vertices of the $n$-dimensional hypercube $\{0,1\}^n$ with $k$-dimensional
subcubes as winning sets.  We describe a pairing strategy which allows \B\ to win if $n$ is a power of 4 and $k \ge n/4 +1$. 
Our results also imply that for all $n \geq 3$  there is a \BPS\  if $k \ge \left\lfloor\frac{3}{7}n\right\rfloor +1$.
\end{abstract}

\section{Introduction}
\label{intro}

A positional game can be thought of as a generalization of \ttt\ played on a hypergraph $(V,\cH)$ where the vertices can be considered the ``board" on which the game is played, and the
{hyper}edges  can be considered the ``winning sets." 
More formally, a {\it positional game} on $(V,\cH)$ is a two-player game where at every  turn each player 
alternately occupies
a previously unoccupied vertex from $V$.  
In a {\it strong positional game}, the first player to occupy all 
vertices of some hyperedge $A\in\cH$ wins.
If at the end of play no hyperedge is completely occupied by either player, that play is declared a draw.
Normal $3\times 3$ \ttt\ is a strong positional game where the vertices
of the hypergraph are the nine positions and the hyperedges are the eight winning lines. 
In a {\it \MB} positional game, the first player, {\it \M}, wins if \M\  
occupies all vertices of some
hyperedge $A\in\cH$, otherwise the second player, {\it \B}, wins.  Therefore, by definition there are no draw plays in \MB\ games.
We say that a player $P$ has a {\it winning strategy} if no matter how the other player
plays, player $P$ wins by following that strategy. 
It is well-known that in a finite \MB\ game, exactly one player has a winning strategy.
(For a nice introduction to positional games, please see \cite{TicTacToeBook}, \cite{BeckInevRandom},
and \cite{PosGamesBook}.)

Recall that the $n$-dimensional  (Boolean) {\em hypercube} $Q_n$ is a bipartite graph 
whose vertex set is $\{0,1\}^n$ and whose edge set is the set of all pairs of vertices
that differ in exactly one coordinate.
A $k$-dimensional {\em subcube}
of $Q_n$ is a subgraph of $Q_n$ that is isomorphic to $Q_k$.
Let $\cQ(n,k)$ denote the hypergraph whose vertex set is $\{0,1\}^n$ and whose hyperedge set is the set of all 
$k$-dimensional subcubes of $Q_n$
(technically, each hyperedge is the set of vertices of some $k$-dimensional subcube).
 For example, $\cQ(n,1) = Q_n$.

In~\cite{MBhypercube}, Kruczek and \sundberg\  initiated the study of the \MB\ game
played on $\cQ(n,k)$.  Using a general \B's win criterion (Erd\H{o}s--Selfridge theorem~\cite{ref:ES})
and \M's win criterion (Theorem~1.2~\cite{TicTacToeBook}), 
they showed that \B\ has a winning strategy for the $\cQ(n,k)$ game when $k \geq \log_2(n) +1$
and \M\ has a winning strategy when $k \leq \log_2 \log_2 (n) -1$.
They also studied the specific question,  ``when can \B\ win by using 
a pairing strategy?"
A {\it pairing strategy} for  \B\
is a set $P$ of pairwise disjoint pairs of vertices in $Q_n$.
\B\  uses $P$ as a strategy by playing as follows:
each time \M\ occupies a vertex $x$, 
if there is an unoccupied vertex $y$ such that $\{x,y\} \in P$, 
then  \B\ immediately responds by occupying $y$; 
otherwise, \B\ occupies an arbitrary unoccupied vertex.  
This guarantees that \B\  occupies at least one vertex from each pair in $P$.  
If every $k$-dimensional subcube of $Q_n$ contains at least one pair from $P$,
then \B\  wins,
and we say $P$ is a {\it Breaker's win} pairing strategy.
(We note that substantial work has been done on \MB\ games on graphs where players occupy edges.
However, in the $\cQ(n,k)$ game, players occupy vertices.  We are unaware of any papers, other than \cite{MBhypercube},
which study the $\cQ(n,k)$ game.)  The goal of this paper is to improve upon the results in
\cite{MBhypercube} pertaining to \BPSs.
%

Let $\ps(n)$ be the smallest value of $k$ such that \B\ can win the positional game on
$\cQ(n,k)$ by {\em using a pairing strategy}.
Proposition~9 in \cite{CPgames} implies that $\ps(n) > \ln(n)$. 
(Indeed, each pair can block at most 
$\binom{n-1}{k-1}$ subcubes and there are at most $2^{n-1}$ pairs.  
When $k = \lfloor \ln(n) \rfloor$, $\binom{n-1}{k-1}2^{n-1} < \binom{n}{k} 2^{n-k}$, which is the total number
of $k$-dimensional subcubes.)
Kruczek and \sundberg~\cite{MBhypercube} showed that $\ps(n) \leq n-3$.
We improve on their result by proving the following:
\newcounter{temp}
\setcounter{temp}{\value{theorem}}
\setcounter{theorem}{3} 
\begin{theorem}
For each $n \geq 3$, there is a \BPS\ for $\cQ(n, \left\lfloor \frac{3}{7} n \right\rfloor +1)$.
\end{theorem}
\setcounter{theorem}{\value{temp}}
 All pairs of vertices in the \BPSs\ that we construct are edges;
thus, our pairing strategies are matchings in $Q_n$.

The remainder of the paper is organized as follows.
In Section~\ref{Q12Five}, we give some basic definitions and explain 
the main techniques
behind constructing our \BPSs\
through an illustrative example.  
In Section~\ref{MainThms}, we state and prove 
a theorem which uses those techniques and can be used to show 
$\ps(n) \leq n/3+1$
if $n = 6\cdot 4^d$ or $n = 9\cdot 4^d$ for some $d \geq 1$.
In Section~\ref{Rotate}, we enhance the techniques from Section~\ref{MainThms}
to prove that 
$\ps(n) \leq n/4+1$ when
$n$ is a power of 4.
In Section~\ref{SpecValues}, we briefly discuss \BPSs\ for specific values of $n$ and $k$, including
the result that $\ps(n) \leq \left\lfloor\frac{3}{7} n \right\rfloor  +1$ for all $n \geq 3$.
In Section~\ref{ExtraResults}, we prove 
$\ps(n) \leq n/3+1$ when
$n$ is a power of 3.
In Section~\ref{conclusion}, we briefly mention how some of our results 
can be viewed as a variation of $d$-polychromatic edge colorings of $Q_n$.

\section{The Basic Strategy}
\label{Q12Five}
The basic idea of our technique is to ``combine" a \BPS\ for $\cQ(4,2)$ with one
for $\cQ(n,k)$ to create a \BPS\ for $\cQ(4n, b)$, where $b = \max \{4k-3, n+k\}$.


We represent each $k$-dimensional subcube of $Q_n$ by a list of $n$ symbols such that
$k$ of the symbols are \Stars\ ($*$) and each of the remaining $n-k$ symbols is 0 or 1.  
For example,  $(*, 0, 1, *)$ represents the set
$$\{0,1\}\times \{0\} \times \{1\}\times \{0,1\} = \{(0,0,1,0),(0,0,1,1),(1,0,1,0),(1,0,1,1)\}$$
which is the vertex set of a 2-dimensional subcube of $Q_4$.   %
We call 0-dimensional and 1-dimensional subcubes {\em vertices} and {\em edges}, respectively.
We abuse terminology and refer to $(*,0,1,*)$ as a ``vector'' with four ``coordinates.''

The following set of edges (vectors) 
is a \BPS\ for $\cQ(4,2)$ (verifying this is straightforward and left to the reader):  %
\begin{align*}
\PS = \{& (*,0,0,0), (0,*,1,0), (0,0,*,1), (0,1,0,*),  
             \\
            &  (*,1,1,1), (1,*,0,1), (1,1,*,0),  (1,0,1,*)\}. 
\end{align*}

{We make use of the following properties in our proofs.}

\medskip

\noindent{\bf Properties of $PS(4,2)$:} 
\begin{enumerate}
\item \label{prop:blocks_2D}
For each pair of indices $1 \leq i < j \leq 4$ and each ordered pair $(b_i, b_j) \in \{0,1\}^2$, %
there is a 
vector  in $\PS$ with $b_i$ (not a \Star) in coordinate~$i$ and $b_j$ (not a \Star) in coordinate~$j$,
because every 2-dimensional subcube in $\cQ(4,2)$ contains an edge from $\PS$.


\item \label{prop:complement}
For each $j\in [4]$, 
$\PS$ has exactly two vectors with a \Star\ in  coordinate~$j$ and those two vectors are
complements of each other (where we consider 0 and 1 to be complements of each other and a \Star\ to be its own complement).  

\end{enumerate}

Given  a \BPS\  $\bv_n$ for $\cQ(n,k)$,
we combine $\bv_n$, as explained below, 
with each edge in $\PS$ to obtain  a set of edges $\binPS$, which we later prove is a \BPS\ 
 for $\cQ(4n,b)$ with $b = \max\{4k-3, n+k\}$. 
Let $\beven_n$ be the set of even parity vectors from $\{0,1\}^n$,
i.e., the sum of the coordinates of each vector is even.
Let $\bodd_n$ be the set of odd parity vectors from $\{0,1\}^n$.
(When it is clear from the context, we drop the subscript $n$, and just write $\beven$ or $\bodd$.)
We take each edge in $\PS$ and replace each 0, 1, or $*$ with 
an element of $\beven_n$, $\bodd_n$, or $\bv_n$, respectively,
to obtain an edge in $\binPS$.
For example, $(0,1,0,*)\in\PS$ yields the following set of edges in $\binPS$:
\begin{equation}\label{PSset}
(\beven_n \times \bodd_n \times \beven_n \times \bv_n) = \{(\vec{x}_1, \vec{x}_2,\vec{x}_3,\vec{x}_4) :
\vec{x}_1 \in \beven_n, \vec{x}_2 \in \bodd_n, \vec{x}_3 \in \beven_n, \vec{x}_4 \in \bv_n \},
\end{equation}
where $(\vec{x}_1, \vec{x}_2,\vec{x}_3,\vec{x}_4)$
represents an edge in $Q_{4n}$, because $\vec{x}_4$ is an edge in $Q_n$ and $\vec{x}_i$ is a vertex in $Q_n$ for $i\in [3]$.  
We call the set in equation~\eqref{PSset} a {\em \binform.}
Doing this for each element of $\PS$, 
and taking the union of the resulting \binform s,
gives the following set of edges: 
\begin{align}
\binPS =\  &(\bv_n \times \beven_n \times \beven_n \times \beven_n) \cup (\beven_n \times \bv_n \times \bodd_n \times \beven_n) \cup  \nonumber \\  
  &(\bv_n \times \bodd_n \times \bodd_n \times \bodd_n)\cup  (\bodd_n \times \bv_n \times \beven_n \times \bodd_n)\cup 
   \nonumber \\
  &(\beven_n \times \beven_n \times \bv_n \times \bodd_n) \cup (\beven_n \times \bodd_n \times \beven_n \times \bv_n)\cup 
  \nonumber   \\   
   &(\bodd_n \times \bodd_n \times \bv_n \times \beven_n)\cup  (\bodd_n \times \beven_n \times \bodd_n \times \bv_n), 
   \label{def:BinPS}
\end{align}
which we later prove is a \BPS\ for $\cQ(4n,b)$ with $b = \max\{4k-3, n+k\}$.

Each subcube $S$ of $Q_{4n}$ can be represented as a length~$4n$ vector with entries from $\{0,1,*\}$.   
We partition the coordinates of that vector into four bins, 
where bin~$j$ contains coordinates $(j-1)n + 1$ through $jn$.
For each subcube $S$ of $Q_{4n}$ and each $j\in [4]$, let $\rS{j}$
denote $S$ {\em restricted} to bin~$j$, i.e., 
if $S = (x_1, \ldots, x_{4n})$, then
$\rS{j} = (x_{(j-1)n+1}, \ldots, x_{jn})$. 
For an edge $\vec{e}\in\binPS$, the sets $\beven_n$, $\bodd_n$, $\bv_n$ determine
the possible values for $\re{j}$; thus, we call $\beven_n$, $\bodd_n$, $\bv_n$
{\em \binsets}.

We say that a subcube 
$S_1$ {\em \blocks} a subcube $S_2$,
if $S_1 \subseteq S_2$.
For example, 
the edge $\vec{e} = (1,0,1,*)$ \blocks\  the subcube $(*,0,1,*)$ 
since $(1,0,1,*) \subset (*,0,1,*)$.
Moreover, $\vec{e}$ blocks exactly three 2-dimensional subcubes of $Q_4$, namely
$(*,0,1,*), (1,*,1,*),$ and  $(1,0,*,*)$.  
We say that a set of subcubes $T$ {\em \blocks} a subcube $S$, 
if $T$ contains a subcube that \blocks\ $S$.
For example, 
since $Q_n$ is a bipartite graph with partite sets $\beven_n$ and $\bodd_n$,
for each edge $\vec{e}$, there exist vertices $\vec{v}_0 \in \beven_n$ and $\vec{v}_1\in \bodd_n$
such that $\vec{v}_0$ \blocks\ $\vec{e}$ and $\vec{v}_1$ blocks $\vec{e}$. Thus,
$\beven_n$ \blocks\ every 1-dimensional subcube of $Q_n$, as does $\bodd_n$.

\medskip
\noindent{\bf Example~1:}
We can form $\binPS[1]$ by using $\beven_1 = \{0\}$, $\bodd_1 = \{1\}$, $\bv_1 = \{ * \}$, where $*$ is the edge $\{0,1\}$,
thus, $\bv_1$ is a \BPS\ for $\cQ(1,1)$.  In this case, $b = \max \{4(1)-3, 1+1\} = 2$, and 
$\binPS[1]$ is precisely $\PS$.

\medskip
\noindent{\bf Example~2:} We can form $\binPS[2]$ by using $\beven_2 = \{(0,0), (1,1) \}$, $\bodd_2 = \{(0,1), (1,0)\}$, 
$\bv_2 = \{ (*,0) \}$, where $\bv_2$ is a \BPS\ for $\cQ(2,2)$.  In this case, $b = \max \{4(2)-3, 2+2 \} = 5$.
{Theorem~\ref{4xThm} implies that $\binPS[2]$ is a \BPS\ for $\cQ(8,5)$. (We directly justify this fact below.)
}

Each \binform\ has cardinality~8, for example,
\begin{align*}
(\beven_2 \times \bodd_2 \times \beven_2 \times \bv_2) =  \{ &(0, 0,\ 0,1, \ 0, 0,\ *, 0 ) ,\   (1,1,\ 0,1, \ 0, 0,\ *, 0 ) \\
    &(0, 0,\ 0,1, \ 1, 1,\ *, 0 ) ,\    (1,1,\ 0,1, \ 1, 1,\ *, 0 ) \\
    &(0, 0,\ 1,0, \ 0,0 ,\ *, 0 ) ,\    (1,1,\ 1,0, \ 0,0 ,\ *, 0 ) \\   
    &(0, 0,\ 1,0, \ 1,1 ,\ *, 0 ) , \   (1,1,\ 1,0, \ 1,1 ,\ *, 0 )\}, 
\end{align*}
(we include extra spaces in the vectors to highlight the four bins)
and $|\binPS[2]| = 64$.

Let us show that $\binPS[2]$ is a \BPS\ for $\cQ(8,5)$.
Let $S$ be a 5-dimensional subcube of $Q_8$.
Recall that $\beven_2$ and $\bodd_2$ each block any subcube of $Q_2$ with positive dimension,
and $\bv_2$ blocks the 2-dimensional subcube of $Q_2$.

{\em Case~1:}  For each $j \in [4]$, the dimension of $\rS{j}$ is positive.
Then, for some $j_1, j_2, j_3, j_4$ with $\{j _1, j_2, j_3, j_4 \} = [4]$,
 $\rS{j_i}$ has dimension~1 for $i\in [3]$ and $\rS{j_4}$ has dimension~2.
We select a \binform\ that has $\bv_2$ in bin~$j_4$ 
{because  $\bv_2$ is only guaranteed to \block\ the 2-dimensional subcube of $Q_2$.}
In each of the other bins, the \binform\ has $\beven_2$ or $\bodd_2$.  
This \binform\ blocks $S$.
For example, if $S = (*, 1, \ 0, *,\ *, *, \ 0, *)$, then $j_4 =3$,
and we can use  $(\beven_2 \times \beven_2 \times \bv_2 \times \bodd_2)$
to block $S$.
Indeed, $(1,1, \ 0,0, \ *,0, \ 0,1) \in (\beven_2 \times \beven_2 \times \bv_2 \times \bodd_2)$
and $(1,1, \ 0,0, \ *,0, \ 0,1) \subset S$.
We could also use $( \bodd_2 \times \bodd_2  \times \bv_2 \times \beven_2)$ to block $S$.

{\em Case~2:}   Suppose $\rS{j_1}$ has dimension~0 for some $j_1\in [4]$.
Then for some $j_2, j_3, j_4$ such that $\{j _1, j_2, j_3, j_4 \} = [4]$,
$\rS{j_2}$ has dimension~1, while $\rS{j_3}$ and $\rS{j_4}$ both have dimension~2.  
{
Because $\bv_2$ is only guaranteed to \block\ the 2-dimensional subcube of $Q_2$,
we must select a \binform\ that has $\bv_2$ in bin~$j_3$ or $j_4$.
W.l.o.g., we will select a \binform\ with $\bv_2$ in bin~$j_4$.
By Property~\ref{prop:complement} of $\PS$, there is a \binform\ $A$
that has $\bv_2$ in bin~$j_4$ and
has a \binset\ ($\beven_2$ or $\bodd_2$) 
whose parity equals the parity of $\rS{j_1}$ in bin~$j_1$. 
We claim that $A$ \blocks\ $S$.}
Bins~$j_2$ and $j_3$ 
{of $A$} 
each contain either $\beven_2$ or $\bodd_2$,
{thus, $A$ \blocks\ $S$.}
For example, if $S = (*, *,\ *, *,\ 1,0, \ *, 1)$, then $j_1 = 3$ and $j_4 \in \{1, 2\}$.
Thus, we use the \binform\ $(\bv_2 \times \bodd_2 \times \bodd_2 \times \bodd_2 )$
to block $S$ because it has $\bodd_2$ in bin~3 and $\bv_2$ in bin~1, or we could
use $(\beven_2 \times \bv_2 \times \bodd_2 \times \beven_2)$
since it has $\bodd_2$ in bin~3 and $\bv_2$ in bin~2.

\medskip
\noindent{\bf Example~3:} We can form $\binPS[3]$ by using $\beven_3 = \{(0,0,0), (0,1,1), (1,0,1), (1,1,0) \}$, 
$\bodd_3 = \{(0,0,1), (0,1,0), (1, 0,0), (1, 1, 1)\}$, 
$\bv_3 = \{ (*,0,0), (1,*,1), (0,1,*) \}$, where $\bv_3$ is a \BPS\ for $\cQ(3,2)$.  In this case, $b = \max \{4(2)-3, 3+2 \} = 5$.
Theorem~\ref{4xThm} implies that $\binPS[3]$ is a \BPS\ for $\cQ(12,5)$.

\section{Main Theorem}
\label{MainThms}

In this section, we prove that $\binPS$ is a matching and  that it blocks
all subcubes of dimension~$b = \max \{4k-3, n+k\}$.

Observe that two subcubes of $Q_n$ are disjoint if and only if for some $j\in [n]$, one of the subcubes
has 0 and the other has 1 in coordinate~$j$ of their vector representations.
\begin{lemma}\label{binPSmatching}
For each $n \geq 1$, 
the set of edges $\binPS$,
as defined in equation~\eqref{def:BinPS},
is a matching.
\end{lemma}

\noindent{\bf Proof of Lemma~\ref{binPSmatching}:}
Let $\vec{e}_1, \vec{e}_2 \in \binPS$ satisfy $\vec{e}_1 \neq \vec{e}_2$.
Then $\rSet{\vec{e}_1}{j} \neq \rSet{\vec{e}_2}{j}$ for some $j\in [4]$.
Let $(x^{(i)}_1, \ldots, x^{(i)}_{4n})$ be the vector representation
of $\vec{e}_i$ for $i\in \{1,2\}$.

{\em Case~1:} Suppose $\rSet{\vec{e}_1}{j}$ and $\rSet{\vec{e}_2}{j}$ are 0-dimensional.
Since $\rSet{\vec{e}_1}{j} \neq \rSet{\vec{e}_2}{j}$, then clearly
there is some $k\in [n]$ such that $\{x^{(1)}_{(j-1)n+k}, x^{(2)}_{(j-1)n+k} \} = \{0,1\}$,
thus, $\vec{e}_1$ and $\vec{e}_2$ are disjoint.

{\em Case~2:} Suppose $\rSet{\vec{e}_1}{j}$ and $\rSet{\vec{e}_2}{j}$ are 1-dimensional.
Then $\rSet{\vec{e}_1}{j}, \rSet{\vec{e}_2}{j} \in \bv_n$.  Since $\bv_n$
is a matching and $\rSet{\vec{e}_1}{j} \neq \rSet{\vec{e}_2}{j}$, there is 
some $k\in [n]$ such that $\{x^{(1)}_{(j-1)n+k}, x^{(2)}_{(j-1)n+k} \} = \{0,1\}$,
thus, $\vec{e}_1$ and $\vec{e}_2$ are disjoint.

Recall that $\binPS$ is the union of eight pairwise disjoint product sets,
as defined in equation~\eqref{def:BinPS}.
If $\vec{e}_1$ and $\vec{e}_2$ are in the same product-set, then either Case~1 or Case~2 occurs.
If they are in distinct product-sets (which correspond to different edges in $\PS$), 
then because $\PS$ is a matching in $Q_4$, we are guaranteed 
that Case~1 occurs for some (possibly other) $j\in [4]$.
\hfill
$\square$

\begin{theorem}\label{4xThm}
If there exists a \BPS\ for $\cQ(n, k)$,  
then there exists a \BPS\ for 
$\cQ(4n,b)$, where $b = \max \{4k-3, n+k\}$.
\end{theorem}

\noindent{\bf Proof of Theorem~\ref{4xThm}:}
Let $S$ be a subcube of $Q_{4n}$.  Let $i\in [4]$. Suppose that $\rS{i}$ has dimension~$c$.
If $ c \geq k$, then we say that $\rS{i}$ 
is {\em \hd\ (\HD)}. %

We use $\binPS$ as our pairing strategy, where $\beven = \beven_n$, $\bodd = \bodd_n$, and $\bv$ is a 
\BPS\ for $\cQ(n,k)$.
Recall that $\beven$ and $\bodd$ can each \block\
any subcube of $Q_n$ with positive dimension, and
$\bv$ can \block\ any  \HD\ subcube of $Q_n$.  

Let $S$ be a $b$-dimensional 
subcube of $Q_{4n}$.
Since $b \geq 4k-3 > 4(k-1)$,  
$\rS{j}$ is \HD\ for some $j\in [4]$.  
Since  $b \geq n + k$ 
and $k \geq 1$, $\rS{j}$ has  dimension~0
for at most two values of $j\in [4]$.  

{\em Case~1:}  Suppose that $\rS{j}$  has positive dimension for each $j\in [4]$, 
and that $\rS{i}$ is \HD\ for some $i \in [4]$.  
Let $A$  be a \binform\ with $\bv$ in bin~$i$.
Since $\rS{i}$ is  \HD, 
$\bv$  \blocks\ $\rS{i}$. 
In each of the other bins, $A$ has either  $\beven$
or  $\bodd$.  
Since $\rS{j}$ has positive dimension for all $j\neq i$, 
 $\beven$ and $\bodd$ can each \block\ $\rS{j}$
 for all $j\neq i$.  Thus, $A$ \blocks\ $S$.

{\em Case~2:}  Suppose instead that $\rS{i}$ and $\rS{j}$ both have dimension~0
for some pair $\{i,j\} \subset [4]$.  
Because of Property~\ref{prop:blocks_2D} of $\PS$,
there is a \binform\  
$A$ whose bin-sets
in bins~$i$ and $j$ match the parities of $\rS{i}$ and $\rS{j}$, respectively. 
There are $b$ coordinates in the vector representation of $S$ that are \Stars.
Since $b \geq n+k$,   $\rS{\ell}$ is \HD\ for each $\ell \in[4]-\{i,j\}$,
and can be \blocked\ by any bin-set $\beven, \bodd,$ or $\bv$.
Therefore, $A$ \blocks\ $S$.

{\em Case~3:} Suppose that $j$ is the unique value in $\{1,2,3,4\}$ such that $\rS{j}$ has dimension~0.  
Also suppose that $\rS{i}$ is \HD\ for some $i\in [4]$.
By Property~\ref{prop:complement} of $\PS$, we can deduce that 
there is a \binform\ $A$  that has $\bv$ in bin~$i$ and a bin-set in bin~$j$
whose parity matches that of $\rS{j}$.
In each of the other two bins, $A$ has either  $\beven$ or  $\bodd$.
Since  $\rS{\ell}$ has positive dimension for each $\ell \in[4]-\{i,j\}$,
$A$ \blocks\ $S$.

In light of Lemma~\ref{binPSmatching}, we conclude that $\binPS$
is a \BPS\ for $\cQ(4n,b)$.  
\hfill $\square$

\begin{corollary}\label{SpCaseThm} 
Suppose there exists a \BPS\ for $\cQ(n,k)$.
\begin{enumerate}[(a)]
\item If $k \geq n/3 +1$, then there exists a \BPS\ for $\cQ(4n, 4k-3)$.
\item If $k = \lfloor n/3 \rfloor +1$, then there exists a \BPS\ for $\cQ(4n, \lfloor 4n/3 \rfloor +1)$.
\end{enumerate}
\end{corollary}

\section{Rotating the Pairing Strategies}
\label{Rotate}

Under certain conditions, we are able to improve upon Theorem~\ref{4xThm} using a technique 
we call, ``rotating pairing strategies."

To illustrate the idea, we construct a \BPS\ for $\cQ(9,4)$ using
the rotating pairing strategy idea.  (The reader may choose to skip to 
Theorem~\ref{RotatePS}.)  For simplicity, we have chosen
an example  based on $\bv_3 = \{ (*,0,0), (1,*,1), (0,1,*) \}$ (instead of $\PS$) 
which has 
only three bins, each of size three.  
We obtain four different \BPSs\ for $\cQ(3,2)$ which together partition the set of edges of $Q_3$,
by using translations of $\bv_3$, namely, 
\begin{align*}
\bvP{0} &= \{(*,0,0), (1,*,1), (0,1,*) \} = \bv_3, \\
\bvP{1} &= \{(*,0,1), (0,*,0), (1,1,*) \} = \bv_3 + (1,0,1),  \\
\bvP{2} &= \{(*,1,0), (0,*,1), (1,0,*) \} = \bv_3 + (1,1,0),  \\
\bvP{3} &= \{(*,1,1), (1,*,0), (0,0,*) \} = \bv_3 + (0,1,1).
\end{align*}
(We will use the fact that  $\bva\cup\bvb\cup\bvc\cup\bvd$ contains every edge in $Q_3$.) 
We partition $\beven_3$ and $\bodd_3$ into four sets each as follows:
\begin{align*}
\bevenP{0} &= \{(0,0,0)\},  & \boddP{0} &= \{(0,0,1) \}, \\
\bevenP{1} &= \{(0,1,1)\},  & \boddP{1} &= \{(0,1,0) \}, \\
\bevenP{2} &= \{(1,0,1)\},  & \boddP{2} &= \{(1,0,0) \}, \\
\bevenP{3} &= \{(1,1,0)\},  & \boddP{3} &= \{(1,1,1) \}. \\
\end{align*}
(In higher dimensions, $\bevenP{j}$ and $\boddP{j}$ will not be singletons.)
Our pairing strategy is 
$$\binPSR[3][3] = (\bv \times \beven \times \beven)^R \cup (\bodd \times \bv \times \bodd)^R \cup (\beven \times \bodd \times \bv)^R, $$
where we define the ``rotating" \binforms\
as follows:
\begin{align*}
(\bv \times \beven \times \beven)^R = \bigcup_{i, j}\ (\bvP{i+j} \times \bevenP{i} \times \bevenP{j}), \\
(\bodd \times \bv \times \bodd)^R = \bigcup_{i, j}\ (\boddP{j} \times \bvP{i+j} \times \boddP{i}), \\
(\beven \times \bodd \times \bv )^R = \bigcup_{i, j}\ (\bevenP{i} \times \boddP{j} \times \bvP{i+j}), 
\end{align*}
where $i, j,$ and $i+j$ are all evaluated modulo~4.  
Notice, for example, that $(\beven \times \bodd \times \bv)^R$ is the union of $4^2$ sets, such as, 
$$( \bevenP{2} \times \boddP{3} \times \bvP{1}) = \{(1,0,1)\} \times \{(1,1,1) \} \times \{(*,0,1), (0,*,0), (1,1,*) \}$$
and 
$$( \bevenP{3} \times \boddP{1} \times \bvP{0}) = \{(1,1,0)\} \times \{(0,1,0) \} \times \{(*,0,0), (1,*,1), (0,1,*) \}.$$

We present an equivalent description of the ``rotating" \binforms.
For $\vec{x}\in \bevenP{j}$ or $\vec{x}\in \boddP{j}$, let $\Index(\vec{x}) = j$.
Then we can write the ``rotating" \binforms\ as 
\begin{align*}
(\bv \times \beven \times \beven)^R &=\{(\vec{x}, \vec{y}, \vec{z}) : \vec{x} \in \bvj, \vec{y}\in \beven, \vec{z}\in\beven,
\mbox{where $j = \Index(\vec{y})+\Index(\vec{z}) \pmod{4}$}\},         \\ 
(\bodd \times \bv \times \bodd)^R &= \{(\vec{x}, \vec{y}, \vec{z}) : \vec{x} \in \bodd, \vec{y}\in \bvj, \vec{z}\in\bodd,
\mbox{where $j = \Index(\vec{x})+\Index(\vec{z}) \pmod{4}$}\},         \\
(\beven \times \bodd \times \bv )^R &= \{(\vec{x}, \vec{y}, \vec{z}) : \vec{x} \in \beven, \vec{y}\in \bodd, \vec{z}\in\bvj,
\mbox{where $j = \Index(\vec{x})+\Index(\vec{y}) \pmod{4}$}\}.
\end{align*}

Theorem~\ref{Q3d_d+1} implies that $\binPSR[3][3]$ is a \BPS\ for $\cQ(9,4)$.
In contrast, 
$$\binPS[3][3] = (\bv_3 \times \beven_3 \times \beven_3) \cup (\bodd_3 \times \bv_3 \times \bodd_3) 
\cup (\beven_3 \times \bodd_3 \times \bv_3),$$
is a \BPS\ for $\cQ(9,5)$ (proof omitted), but not for $\cQ(9,4)$. 
For example, let $S$ be a 4-dimensional subcube of $Q_9$ so that
$\rS{1}$ is 0-dimensional, $\rS{2}$ is 1-dimensional, and $\rS{3}$ is 3-dimensional.
If $\rS{1}$ has even parity, then both $\binPS[3][3]$ and $\binPSR[3][3]$ \block\ $S$ using
the \binforms\ $(\beven_3 \times \bodd_3 \times \bv_3)$ and $(\beven \times \bodd \times \bv )^R$,
respectively.  (This is because $\bodd_3= \bodd$ \blocks\ subcubes with positive dimension,
and any subcube of $Q_3$ \blocks\ $\rS{3}$ because it is 3-dimensional.) 
However, if $\rS{1}$ has odd parity, then $(\bodd \times \bv \times \bodd)^R$ \blocks\ $S$, 
but $(\bodd_3 \times \bv_3 \times \bodd_3)$ might not (because $\rS{2}$ is 1-dimensional, but bin~2 contains $\bv_3$).
For example, suppose
$S = (0,1,0, \ *,1,0, \ *,*,*)$. 
Observe that $(\boddP{1} \times  \bvP{2} \times  \boddP{1}) \subset (\bodd \times \bv \times \bodd)^R$,
and $(0,1,0, \ *,1,0,\ 0,1,0) \in (\boddP{1} \times  \bvP{2} \times  \boddP{1})$, which \blocks\ $S$.
However, since $\rS{2} = (*,1,0)$, which is not \block ed by $\bv_3= \{ (*,0,0), (1,*,1), (0,1,*) \}$,
$(\bodd_3 \times \bv_3 \times \bodd_3)$ does not \block\ $S$.  

\begin{theorem}\label{RotatePS}
Suppose there exists a set of matchings $\{\bva, \ldots, \bv^{(m-1)}\}$ 
such that each $\bvj$ is a 
\BPS\ for $\cQ(n,k)$ 
and $\bigcup_j \bvj$ equals the set of edges of $Q_n$.
Moreover, suppose that 
there is a partition of $\beven_{n}$ (of $\bodd_n$)
of size~$m$
such that every subcube of $Q_n$ of dimension $n-k+2$ contains at least 
one vertex from each of the sets in the partition of $\beven_n$ (of $\bodd_n$).
Then there exists a \BPS\ for $\cQ(4n, b)$, where $b = \max \{4k-3, n+1\}$.
\end{theorem}

\noindent{\bf Proof of Theorem~\ref{RotatePS}:}

Let $\{\bevenP{0}, \ldots \bevenP{m-1}\}$ and $\{\boddP{0}, \ldots, \boddP{m-1}\}$
be the given partitions of $\beven_n$ and $\bodd_n$, respectively.
For every $\vec{x}\in \bevenP{j} \cup  \boddP{j}$, we define $\Index(\vec{x})$ to be $j$.
Let 
\begin{align*}
\binPSR = \\
& (\bv \times \beven \times  \beven \times  \beven)^R \cup  (\beven \times  \bv \times  \bodd \times  \beven)^R\ \cup  \\
& (\bv \times \bodd \times  \bodd \times  \bodd)^R \cup   (\bodd \times  \bv \times  \beven \times  \bodd)^R\ \cup  \\
& (\beven \times  \beven \times  \bv \times  \bodd)^R \cup  (\beven \times  \bodd \times \beven \times  \bv)^R\ \cup    \\ 
& (\bodd \times  \bodd \times  \bv \times  \beven)^R \cup   (\bodd \times  \beven \times \bodd \times  \bv)^R
\end{align*}
where, for example,
\begin{align*}
(\beven \times  \beven \times  \bv \times  \bodd)^R = 
\{(\vec{x}, \vec{y}, \vec{z}, \vec{w}) :\  & \vec{x} \in \beven, \vec{y}\in \beven, \vec{z}\in\bvj, \vec{w} \in \bodd, \\
& \mbox{where } j = \Index(\vec{x})+\Index(\vec{y}) + \Index(\vec{w}) \;  \pmod{m} \}.
\end{align*}

Let $b = \max\{4k-3,n+1\}$, and let $S$ be a $b$-dimensional subcube of $Q_{4n}$.
Suppose that $\rS{i}$ has dimension~$c$ for some $i\in [4]$.
If $1 \leq c \leq k-1$, then 
we say that $\rS{i}$ 
is {\em \ld\ (\LD)}. 
(Note: if $\rS{i}$ has dimension~0, it is {\em not} LD.)
Recall, if $ c \geq k$, then  $\rS{i}$ 
is \HD. 

Since $b \geq 4k-3$, $\rS{i}$ is \HD\ for at least one $i \in [4]$. 
Since $b \geq n+1$, $\rS{i}$ has dimension~0 for at most two values of $i\in [4]$. 

We can follow the proof of Theorem~\ref{4xThm} for Case~1 and Case~3 in that proof.
For Case~2 in that proof, if 
$\rS{i_1}$ and $\rS{i_2}$ both have dimension~0, and $\rS{i_3}$ and $\rS{i_4}$
are both \HD, where $\{i_1, i_2, i_3, i_4\} = [4]$, then we can follow the proof of Theorem~\ref{4xThm}.
Therefore, we may assume that $\rS{i_1}$ and $\rS{i_2}$ both have dimension~0,
 $\rS{i_3}$ is \LD, and $\rS{i_4}$ is \HD.
Property~\ref{prop:blocks_2D} of $\PS$ implies that there is a \binform\ $A$ whose bin-sets
in bins~$i_1$ and $i_2$ match the parities of $\rS{i_1}$ and $\rS{i_2}$, respectively.
If  
$A$ has  $\beven$ or $\bodd$ in  bin~$i_3$,  
then $A$ \blocks\ $S$ because $\beven$ and $\bodd$ each  \block\ 
subcubes with any positive dimension, and $\bv$ \blocks\ \HD\ subcubes.
Suppose instead that $A$ has $\bv$ in bin~$i_3$,
and w.l.o.g., $A$ has $\beven$ in bin~$i_4$. 
Since $\bigcup_j \bvj$ equals the set of edges of $Q_n$
and $\rS{i_3}$ has positive dimension, 
there is a bin-set $\bvj$ which contains an edge that  \blocks\ $\rS{i_3}$.  
Let $c\in \{0,\ldots, m-1\}$ satisfy $\Index(\rS{i_1}) + \Index(\rS{i_2}) + c  = j~\pmod{m}$.
Since, by assumption, 
every subcube of $Q_n$ of dimension $n-k+2$ contains at least one vertex from $\bevenP{c}$,
there is a vertex $\vec{x} \in \bevenP{c}$ that \blocks\ 
$\rS{i_4}$ as long as $\rS{i_4}$ has dimension at least $n-k+2$.   
Since $\rS{i_3}$ is \LD\ 
and the sum of the dimensions of $\rS{i_3}$ and $\rS{i_4}$ is $b \ge n+1$, 
$\rS{i_4}$ has dimension at least $n+1-(k-1) = n-k+2$. 
Therefore, we can find an edge in $A$ that \blocks\ $S$.
\hfill $\square$

For Theorem~\ref{RotatePS} to be useful, we need there to exist,
for some $\cQ(n,k)$ with $k < n/3 + 1$, 
a set of \BPSs\ whose union is the set of edges of $Q_n$;
for if $k \geq n/3+1$, then $4k-3 \geq n+k$, and hence $b= 4k-3$ in both
Theorems~\ref{4xThm} and \ref{RotatePS}.
The existence of such sets of pairing strategies  is proved in Theorem~\ref{Q4d_d+1},
which states that, for $d \ge 0$, the edges of $Q_{4^{d+1}}$ can be partitioned
into $4^{d+1}$ \BPSs\ for $\cQ(4^{d+1}, 4^d +1)$.
We prove Theorem~\ref{Q4d_d+1} by induction on $d$.
To introduce definitions and techniques we use in that proof,
let us consider here the cases $d=0$ and $d=1$. 

We have seen that $\PS$ is a \BPS\ for the case $d=0$.
To  obtain four different \BPSs\ for $\cQ(4,2)$ which together partition the set of edges of $Q_4$,
we use translations of $\PS$, namely,
\begin{align*}
PS_0(4,2) &= \PS &=  & & \{ &(*,0,0,0), (0,*,1,0),  (0,0,*,1), (0,1,0,*),  
             \\
   &    &  &  &   &  (*,1,1,1),  (1,*,0,1), (1,1,*,0), (1,0,1,*)\}, 
            \\ \\
PS_1(4,2) &= \PS + (0,0,1,1) &=&  & \{& (*,0,1,1),  (0,*,0,1), (0,0,*,0), (0,1,1,*),  
             \\
   &  &   &  &    &  (*,1,0,0), (1,*,1,0), (1,1,*,1), (1,0,0,*)\}, 
            \\ \\
PS_2(4,2) &= \PS + (0,1,0,1) &= & &\{  &(*,1,0,1),  (0,*,1,1), (0,1,*,0), (0,0,0,*),  
             \\
 &   &   &  &     &  (*,0,1,0), (1,*,0,0), (1,0,*,1),  (1,1,1,*)\}, 
            \\ \\
PS_3(4,2) &= \PS + (0,1,1,0) &= & & \{ &(*,1,1,0), (0,*,0,0), (0,1,*,1), (0,0,1,*),  
             \\
 &    &  &  &    &  (*,0,0,1), (1,*,1,1), (1,0,*,0), (1,1,0,*)\}. 
\end{align*}
Since every translation is an automorphism,  
$PS_j(4,2)$, for $0 \leq j \leq 3$, is a \BPS\ for $\cQ(4,2)$  that satisfies Properties~1 and 2 of $\PS$.
  Thus, we now have four \BPSs\ from which
to construct our rotating \binforms\ for $d>0$.

For the case $d=1$,
we can obtain one \BPS\ for $\cQ(16,5)$ 
by applying Theorem~\ref{RotatePS},  with $\bvj = PS_j(4,2)$ for $0 \leq j \leq 3$ and
using any partitions of $\beven_4$ and $\bodd_4$ with four nonempty parts each (since $n-k+2 = 4$).
However, to obtain sixteen \BPSs\ for
$\cQ(16,5)$ which together partition the set of edges of $Q_{16}$, 
Theorem~\ref{Q4d_d+1} generalizes the idea of Theorem~\ref{RotatePS} by 
using ``shifted" rotating \binforms, which we define below.

Suppose $\bigcup_0^{m-1} \bvP{j}$ equals the set of edges in $Q_n$, 
with each $\bvP{j}$ a \BPS\ for $\cQ(n,k)$,
and $\{\bevenP{0}, \ldots \bevenP{m-1}\}$ and $\{\boddP{0}, \ldots, \boddP{m-1}\}$
are partitions of $\beven_n$ and $\bodd_n$, respectively.
An example of a {\em rotating \binform\ shifted by} $s$, where $0 \leq s \leq m-1$, is
\begin{align*}
(\beven \times  \beven \times  \bv \times  \bodd)^{\Rs} = \{(\vec{x},& \vec{y}, \vec{z}, \vec{w}) :\   \vec{x} \in \beven, \vec{y}\in \beven, 
\vec{z}\in\bvP{i}, \vec{w} \in \bodd, \\
& \mbox{where $i = s + \Index(\vec{x})+\Index(\vec{y}) + \Index(\vec{w}) \; \pmod{m}$}\}.
\end{align*}
Note that when $s=0$, we obtain $(\beven \times  \beven \times  \bv \times  \bodd)^{R(0)} = (\beven \times  \beven \times  \bv \times  \bodd)^{R}$.

For each $s\in \{0,\ldots,m-1\}$, we define four pairing strategies
$\binPSP{0}{s}$, $\binPSP{1}{s}$, $\binPSP{2}{s}$, $\binPSP{3}{s}$ based on 
$PS_0(4,2)$, $PS_1(4,2)$, $PS_2(4,2)$, $PS_3(4,2)$, respectively. 
For example, 
\begin{align*}
&\binPSP{3}{s} = \\
& (\bv \times \bodd \times  \bodd \times  \beven)^{\Rs} \cup  (\beven \times  \bv \times  \beven \times  \beven)^{\Rs} \cup  
(\beven \times   \bodd \times  \bv \times  \bodd)^{\Rs} \cup  (\beven \times  \beven \times \bodd \times  \bv)^{\Rs}    
 \\ 
&\cup  (\bv \times \beven \times  \beven \times  \bodd)^{\Rs} \cup   (\bodd \times  \bv \times  \bodd \times  \bodd)^{\Rs} \cup  
 (\bodd \times  \beven \times  \bv \times  \beven)^{\Rs} \cup   (\bodd \times  \bodd \times \beven \times  \bv)^{\Rs}.
\end{align*}

Thus, the sixteen \BPSs\ for
$\cQ(16,5)$ which result from Theorem~\ref{Q4d_d+1}
are $\binPSP[4]{j}{s}$ for $0 \leq j \leq 3$ and $0 \leq s \leq 3$, 
where we use the bin-sets $\beven_4$, $\bodd_4$ and $\bvP{j} = PS_j(4,2)$
in the definition of each $\binPSP[4]{j}{s}$.


Theorem~\ref{Q4d_d+1} requires the following lemma.

\begin{lemma}\label{HalfPlus1}
For all $k \geq 1$, the sets $\beven_{2^k}$ and $\bodd_{2^k}$ can be partitioned into 
subsets  $A_1, \ldots, A_{2^k}$ 
and $B_1, \ldots, B_{2^k}$, 
respectively,
so that every subcube of $Q_{2^k}$ of dimension $2^{k-1} + 1$
contains a vertex in every $A_j$ and $B_j$. 
\end{lemma}

\noindent{\bf Proof of Lemma~\ref{HalfPlus1}:}
We proceed by induction on $k$.  The case $k=1$ is trivial.
Suppose for some $k \ge 1$ we have 
$A_1, \ldots, A_{2^{k}}$ and $B_1, \ldots, B_{2^{k}}$
as in the statement of the lemma.
For each $\ell  \in [2^{k}]$, let 
$$G_{\ell} = \bigcup_{i + j \equiv \ell} (A_i \times A_j)$$ 
$$H_{\ell} = \bigcup_{i + j \equiv \ell} (B_i \times B_j)$$
$$I_{\ell} = \bigcup_{i + j \equiv \ell} (A_i \times B_j)$$ 
$$J_{\ell} = \bigcup_{i + j \equiv \ell} (B_i \times A_j)$$
where all the equivalencies are taken modulo $2^k$.
We will show that $\{G_1, \ldots, G_{2^{k}},$ $H_1, \ldots, H_{2^{k}} \}$ is the desired partition
of $\beven_{2^{k+1}}$. A similar proof shows that $\{I_1, \ldots, I_{2^{k}},$ $J_1, \ldots, J_{2^{k}} \}$ is the desired partition of $\bodd_{2^{k+1}}$.

We can write any vertex $v \in Q_{2^{k+1}}$
 as $v = (v_1, v_2)$ with $v_i \in Q_{2^k}$.
Let $v \in \beven_{2^{k+1}}$;
then $v_1$ and $v_2$ must have the same parity,
 so $(v_1, v_2)$ is in either $(A_i \times A_j)$ or $(B_i \times B_j)$ for some $i, j \in [2^k]$,
 and hence $v$ is in either $G_{\ell}$ or $H_{\ell}$,  where $i + j \equiv \ell \pmod{2^k}$.
 As these sets are pairwise disjoint,
$\{G_1, \ldots, G_{2^{k}},$ $H_1, \ldots, H_{2^{k}} \}$ partitions $\beven_{2^{k+1}}$.

Given a subcube $S$ of  $Q_{2^{k+1}}$,
divide its coordinates into two bins, 
with  coordinates  $1, \ldots, 2^{k}$ in bin~1, 
and $2^{k} +1, \ldots, 2^{k+1}$ in bin~2.
If $S$ has dimension $2^k +1$,
then $\rS{1}$ or $\rS{2}$, say $\rS{1}$,  has dimension at least $2^{k-1} + 1$.
So $\rS{1}$ contains a vertex in every $A_i$ and $B_i$.
And $\rS{2}$ has dimension at least one,
so it contains a vertex in each of  $\beven_{2^{k}}$ and $\bodd_{2^{k}}$,
i.e., in some $A_j$ and some $B_{j'}$.
Now, for every $\ell \in [2^k]$
there exist $i$ and $i'$ such that
$i+j \equiv i'+j' \equiv \ell \pmod{2^k}$.
So $S$ contains a vertex in each of 
$(A_i \times A_j)$, $(B_i \times A_j)$, $(A_{i'} \times B_{j'})$, and $(B_{i'} \times B_{j'})$.
So $S$ contains a vertex in every $G_{\ell}$ and $H_{\ell}$.
 \hfill $\square$

\begin{theorem}\label{Q4d_d+1}
For each $d \geq 0$, there exist $4^{d+1}$ pairwise disjoint \BPSs\ for $\cQ(4^{d+1}, 4^d +1)$ 
with equal cardinalities
which partition the set of edges of $Q_{4^{d+1}}$.
\end{theorem}

\noindent{\bf Proof of Theorem~\ref{Q4d_d+1}:}
We proceed by induction on $d$.  The \BPSs\
$PS_j(4,2)$ for $0 \leq j \leq 3$ handle the case $d=0$.
Let $d \geq 1$.  By the inductive hypothesis,
there exist $4^d$ disjoint \BPSs\ $\bvP{0}, \ldots, \bvP{4^d -1}$ for
$\cQ(n,k) = \cQ(4^{d}, 4^{d-1} +1)$ with equal cardinalities 
which partition the set of edges of $Q_{4^{d}}$.
We will show that $\binPSP{j}{s}$ is a \BPS\
for $\cQ(4^{d+1}, 4^d +1)$ for $0 \leq j \leq 3$ and $0 \leq s \leq 4^d -1$,
where we use $\beven_{4^d}$, $\bodd_{4^d}$, and the \BPSs\ $\bvP{i}$
from the inductive hypothesis in the definitions of the \binform s.
Moreover, we will show that the \BPSs\ $\binPSP{j}{s}$ form a partition
of the edges of $Q_{4^{d+1}}$.
To do so, we follow the proof of Theorem~\ref{RotatePS}, with $m = n = 4^d$,
and with one minor change, as follows.

By Lemma~\ref{HalfPlus1}, 
there is a partition of $\beven_{4^d}$ (of $\bodd_{4^d}$)
of size $m$ such that every subcube of $Q_n$ of dimension 
$\frac{1}{2} 4^d +1$ 
contains at least 
one vertex from each of the sets in the partition of $\beven_{4^d}$ (of $\bodd_{4^d}$).
Since $n - k + 2 = \frac{3}{4} 4^d +1 > \frac{1}{2} 4^d +1$, the hypotheses for 
Theorem~\ref{RotatePS} are satisfied.  
We can substitute $\binPSP{j}{s}$ for $\binPSR$ 
throughout the proof of Theorem~\ref{RotatePS} and reach
the conclusion that $\binPSP{j}{s}$ is a \BPS\ for 
$\cQ(4n, n+1) = \cQ(4^{d+1}, 4^d +1)$ (since $\max(4k-3, n+1)=  n+1 = 4k-3)$.
The only minor change we make is to write
``let $c\in \{0,\ldots, m-1\}$ satisfy $(s + \Index(\rS{i_1}) + \Index(\rS{i_2}) + c) = j \pmod{m} $."

It remains to show that the sets $\binPSP{j}{s}$ form a partition of $E(Q_{4^{d+1}})$.
Let $E(Q_{4^{d+1}})$ be the set of edges of $Q_{4^{d+1}}$.
We will first show that $$E(Q_{4^{d+1}}) \subseteq \bigcup_{j,s}\binPSP{j}{s},$$
which implies $\bigcup_{j,s}\binPSP{j}{s} = E(Q_{4^{d+1}})$.

Let $S \in E(Q_{4^{d+1}})$.   Suppose that $\rS{i_1}$, $\rS{i_2}$, $\rS{i_3}$ are all vertices in $Q_{4^d}$, 
and $\rS{i_4}$ is an edge in $E(Q_{4^d})$. 
Let $\vec{x}$ be 
the edge in $E(Q_4)$ such that
coordinate $i_\ell$ of $\vec{x}$ matches the parity of $\rS{i_\ell}$ for $1 \leq \ell \leq 3$,
and coordinate $i_4$ of $\vec{x}$ is a \Star.  
Since $PS_0(4,2)$, $PS_1(4,2)$, $PS_2(4,2)$, $PS_3(4,2)$ partition $E(Q_4)$, 
$\vec{x} \in PS_j(4,2)$ for some $0 \leq j \leq 3$.

Since $\bvP{0}, \ldots, \bvP{4^d-1}$ partition $E(Q_{4^d})$, 
$\rS{i_4} \in \bvP{t}$ for some $0 \leq t \leq 4^d-1$.
Let $s\in \{0,\ldots, 4^d-1\}$ satisfy
$$s + \Index(\rS{i_1}) + \Index(\rS{i_2})+ \Index(\rS{i_3})= t \pmod{4^d}. $$
Let $A_s$ be the \binform\ from $\binPSP{j}{s}$ which corresponds to $\vec{x}$.
Then $S \in A_s$ because of how $\vec{x}$ was chosen.
Therefore,  $E(Q_{4^{d+1}}) \subseteq \bigcup_{j,s} \binPSP{j}{s}$, as desired.

It remains to show that the sets $\binPSP{j}{s}$ with  $0 \leq j \leq 3$ and  $0 \leq s \leq 4^d -1$ are pairwise disjoint.
Since every $\bvP{i}$ has the same cardinality and 
$\bvP{0}, \ldots, \bvP{4^d-1}$ partition
$E(Q_{4^d})$, which has
cardinality $4^d (2^{4^d -1})$, $|\bvP{i}| = 2^{4^d -1}$ for $0 \leq i \leq 4^d-1$.
Each $\binPSP{j}{s}$ is the union of eight \binform s.  
Since $|\beven| = |\bodd| = 2^{4^d-1}$,  each \binform\ has cardinality
$\left(2^{4^d-1}\right)^3 (2^{4^d -1})$.
Thus, $\left| \binPSP{j}{s}\right| \leq 8 \left(2^{4^d -1}\right)^4 = 2^{4^{d+1} -1},$ and
\begin{equation*}
\left|\bigcup_{j,s} \binPSP{j}{s} \right| \leq \sum_{j,s} \left| \binPSP{j}{s}\right| \leq 4^{d+1} 2^{4^{d+1} -1} 
= \left|E(Q_{4^{d+1}})\right|.
\end{equation*}
Therefore the inequalities must in fact be equalities,
and hence the \BPSs\ $\binPSP{j}{s}$ form a partition of $E(Q_{4^{d+1}})$.
\hfill $\square$

\section{Pairing Strategies for Specific Values of $n$ and $k$}   
\label{SpecValues}

Both Theorems~\ref{4xThm} and \ref{RotatePS} require the existence
of a \BPS\ for a game played on the vertices of $Q_n$ to construct a \BPS\ 
for a game played on the vertices of $Q_{4n}$.  
The following two lemmas allow us to construct \BPSs\ for 
games played on the vertices of $Q_d$ where $d$ is not divisible by 4.  

\begin{lemma}[\cite{MBhypercube}]\label{plus1plus1}
If there is a \BPS\ for the \MB\ game played on $\cQ(n,k)$, then there is a \BPS\
for the \MB\ game played on $\cQ(n+1, k+1)$.
\end{lemma}

\begin{lemma}\label{smallerN}
If there exists a matching which is a \BPS\ for the \MB\ game played on $\cQ(N,k)$, then there is 
a matching which is a \BPS\
for the \MB\ game played on $\cQ(n, k)$ for all $n \leq N$.
\end{lemma}

Both lemmas are fairly easy to justify.  For a full proof of Lemma~\ref{plus1plus1}, see
\cite{MBhypercube}.  To understand the idea behind Lemma~\ref{smallerN}, for example,
observe that 
there is a natural correspondence between the set of $k$-dimensional subcubes of $Q_n$ 
and the set of $k$-dimensional subcubes of $Q_N$ whose last $N-n$ coordinates
are fixed at 0. The set of edges from our \BPS\ which \blocks\ those $k$-dimensional
subcubes must also have their last $N-n$ coordinates fixed at 0.
If we truncate each of those edges after their $n^{th}$ coordinate, we will
obtain a \BPS\ for the set of $k$-dimensional subcubes of $Q_n$.

So far we have exhibited \BPSs\ for 
$\cQ(3,2)$, $\cQ(4,2)$ and $\cQ(9,4)$.  
Let us provide a \BPS\ for $\cQ(6,3)$
in order to help us  construct \BPSs\ for other values of $n$ and $k$.

To construct a \BPS\ for $\cQ(6,3)$, 
we will use sets of edges that
{\em resemble} cyclic permutations.
For example, let
\begin{align*}
\langle (*,0,1,0,0,0)\rangle = \{ &(*,0,1,0,0,0), \\ & (0,*,0,1,0,0),  \\ & (0,0,*,0,1,0), \\ & (0,0,0,*,0,1), 
\\ & (1,0,0,0,*,0), \\ & (0,1,0,0,0,*)\}.
\end{align*}
Then, 
\begin{equation}
\langle (*, 0, 1, 0, 0, 0)\rangle  \cup \langle (*, 1, 0, 1, 1,1)\rangle
\cup \langle (*,0,1,1,0,0)\rangle \cup \langle (*,1,0,0,1,1)\rangle
\end{equation}
is a \BPS\ for $\cQ(6,3)$ consisting of 24 edges (verified by computer).

If we start with our \BPS\ for $\cQ(6,3)$ and  repeated apply Corollary~\ref{SpCaseThm}(b),
then we obtain a \BPS\ for $\cQ(6\cdot 4^n, 2\cdot 4^n +1)$ for all $n \geq 0$.
Likewise, if we start with our \BPS\ for $\cQ(9,4)$ and  repeated apply Corollary~\ref{SpCaseThm}(b),
then we obtain a \BPS\ for $\cQ(9\cdot 4^n, 3\cdot 4^n +1)$ for all  $n \geq 0$.
Theorem~\ref{Q4d_d+1} states that there is a \BPS\
for $\cQ(4^{n+1}, 4^n +1)$ for all $n \geq 0$.

For each $n \geq 0$, there remain three intervals for which we have not yet 
described a \BPS, namely, 
for games played on the vertices of $Q_d$
where $4^{n+1} < d < 6 \cdot 4^n$, or $6\cdot 4^n < d < 9 \cdot 4^n$, or 
$9 \cdot 4^n < d < 4^{n+2}$. 
To establish the existence of \BPSs\ for these values of $d$,
we can use Lemmas~\ref{plus1plus1} and \ref{smallerN}. 
We use the same approach for each interval.  Specifically,
for an interval of the form $N_1 < d < N_2$, we have a \BPS\
for $\cQ(N_1, k_1)$ and $\cQ(N_2, k_2)$.  When $d = N_1 + j$ for $1 \leq j \leq 4^n$,
we use Lemma~\ref{plus1plus1} and our \BPS\ for $\cQ(N_1, k_1)$
to obtain a \BPS\ for $\cQ(N_1+ j, k_1 + j)$.  When $d = N_1 + j$ for $4^n + 1 \leq j < N_2$,
we use Lemma~\ref{smallerN} and our \BPS\ for $\cQ(N_2, k_2)$
to obtain a \BPS\ for $\cQ(N_1 + j, k_2)$.

After applying this technique to all three intervals, we obtain \BPSs\ for 
%
\begin{align*}
&\cQ(4^{n+1} + 0\cdot 4^n + j, \, 1\cdot 4^n +1+ j),  & 1 & \leq j \leq 4^n, \\
&\cQ(4^{n+1} + 1\cdot 4^n + j, \, 2\cdot 4^n +1),  & 1 & \leq j \leq 4^n, \\
&\cQ(4^{n+1} + 2\cdot 4^n + j, \, 2\cdot 4^n +1+ j),  & 1 & \leq j \leq 4^n, \\
&\cQ(4^{n+1} + 3\cdot 4^n + j, \, 3\cdot 4^n +1),   & 1 & \leq j \leq 2\cdot 4^n, \\
&\cQ(4^{n+1} + 5\cdot 4^n + j, \, 3\cdot 4^n +1+ j), & 1 & \leq j \leq 4^n, \\
&\cQ(4^{n+1} + 6\cdot 4^n + j, \, 4\cdot 4^n +1) & 1 & \leq j \leq 6\cdot 4^n.
\end{align*}

Using the results stated in this section, we have established the existence of a non-trivial 
\BPS\ for $\cQ(N,K)$ for each $N \geq 3.$
When $N = 4^{n+1}$, we have that $K$ is  $N/4 +1$.
When $N = 6\cdot 4^n$ or $N = 9 \cdot 4^n$, we have that $K$ is $N/3 +1$.
We can ask the following question. What is the largest value that  the ratio $K/N$ attains using the results above?
We observe that 
as $N$ increases from $4^{n+1}$ to $5\cdot 4^n$, 
the ratio $K/N$ increases from $\frac{1}{4} + \frac{1}{N}$ to $\frac{2}{5} + \frac{1}{N}$.
As $N$ increases from $5\cdot 4^n$ to $6\cdot 4^n$, 
$K/N$ decreases from $\frac{2}{5} + \frac{1}{N}$ to $\frac{1}{3} + \frac{1}{N}$.
As $N$ increases from $6\cdot 4^{n}$ to $7\cdot 4^n$, 
$K/N$ increases from $\frac{1}{3} + \frac{1}{N}$ to $\frac{3}{7} + \frac{1}{N}$.
As $N$ increases from $7\cdot 4^n$ to $9\cdot 4^n$, 
$K/N$ decreases from $\frac{3}{7} + \frac{1}{N}$ to $\frac{1}{3} + \frac{1}{N}$.
As $N$ increases from $9\cdot 4^{n}$ to $10\cdot 4^n$, 
$K/N$ increases from $\frac{1}{3} + \frac{1}{N}$ to $\frac{2}{5} + \frac{1}{N}$.
As $N$ increases from $10\cdot 4^n$ to $4^{n+2}$, 
$K/N$ decreases from $\frac{2}{5} + \frac{1}{N}$ to $\frac{1}{4} + \frac{1}{N}$.
The largest value $K/N$ achieves is $\frac{3}{7} + \frac{1}{N}$, when $N = 7\cdot 4^n$
and $K = 3\cdot 4^n +1$.  One can check that for each $N \geq 3$, there is a 
\BPS\ for $K = \left\lfloor \frac{3}{7} N \right\rfloor + 1$.
Thus, we have the following theorem:
\begin{theorem}
For each $N \geq 3$, there is a \BPS\ for $\cQ(N, \left\lfloor \frac{3}{7} N \right\rfloor +1)$.
\end{theorem}

We present 
the values of $N$ and $K$ corresponding to the (locally) minimum and (locally) maximum values achieved by
$K/N$ in the following table. 
$$
\begin{array}{| c | c | c|}
\hline
N & K & \mbox{(local) max/min} \\ \hline    \hline
4^n & N/4 + 1 & \mbox{min} \\  \hline
5 \cdot 4^n & (2/5)N + 1 & \mbox{max} \\  \hline
6\cdot 4^n & N/3 +1 & \mbox{min} \\   \hline
7 \cdot 4^n & (3/7)N + 1 & \mbox{max} \\  \hline
9\cdot 4^n & N/3 +1 & \mbox{min} \\  \hline
10 \cdot 4^n & (2/5)N + 1 & \mbox{max} \\  \hline
\end{array}
$$


\section{Extra Results}   
\label{ExtraResults}

In Lemma~\ref{PartitionEvenOdd}, we state a generalization of Lemma~\ref{HalfPlus1}.  
Lemma~\ref{PartitionEvenOdd} is used to prove Theorem~\ref{Q3d_d+1}, but it is also interesting in its own right.  
Theorem~\ref{Q3d_d+1} provides additional pairing strategies that are not covered in Section~\ref{SpecValues}.

\begin{lemma}\label{PartitionEvenOdd}
For all $n \geq 1$ and all $c \geq 2$, the sets $\beven_{c^n}$ and $\bodd_{c^n}$ can each be partitioned into 
subsets $A_1, \ldots, A_{(2^{c-1})^n}$ 
and $B_1, \ldots, B_{(2^{c-1})^n}$,  
respectively,
so that every subcube $S$ of $Q_{c^n}$ of dimension $c^n - c^{n-1} + 1$
contains a vertex 
from each $A_j$ and each $B_j$.
\end{lemma}

\noindent{\bf Proof of Lemma~\ref{PartitionEvenOdd}:}
We proceed by induction on $n$.  The case $n=1$ is trivial.
Suppose for some $n \geq 1$ we have $A_1, \ldots, A_{(2^{c-1})^n}$ 
and $B_1, \ldots, B_{(2^{c-1})^n}$ as in the statement of the lemma.

For each vector $(b_1, \ldots, b_c) \in \{0,1\}^c$
and each vector of indices $(i_1, \ldots, i_c) \in [(2^{c-1})^{n}]^c$, 
define the set
$(D_{i_1}\times \cdots \times D_{i_c})$ where $D_{i_j} = A_{i_j}$ if $b_j = 0$
and $D_{i_j} = B_{i_j}$ if $b_j = 1$.  For example, if $c=3$ and $n=1$, we
could have $(1, 1, 0) \in \{0,1\}^3$ and $(4,1,2) \in [4]^3$ which
result in the set $(B_4 \times B_1 \times A_2)$.
Then for each $\vec{b} \in \beven_c$ and each $\ell \in [(2^{c-1})^{n}]$,
we define the set
$$A_{(\vec{b},\ell)} = \bigcup_{i_1+ \cdots + i_c \equiv \ell} (D_{i_1}\times \cdots \times D_{i_c}),$$
where the congruence is taken modulo $(2^{c-1})^{n}$.
For example, in the case $c = 3$, we define $4^{n+1}$ sets
\begin{align*}
A_{((0,0,0),\ell)} &= \bigcup_{i_1+ i_2+ i_3 \equiv \ell} (A_{i_1} \times A_{i_2}\times A_{i_3}), & 1 \leq \ell \leq 4^{n}, \\
A_{((0,1,1),\ell)} &= \bigcup_{i_1+ i_2+ i_3 \equiv \ell} (A_{i_1} \times B_{i_2}\times B_{i_3}), & 1 \leq \ell \leq 4^{n}, \\
A_{((1,0,1),\ell)} &= \bigcup_{i_1+ i_2+ i_3 \equiv \ell} (B_{i_1}\times A_{i_2}\times B_{i_3}), & 1 \leq \ell \leq 4^{n}, \\
A_{((1,1,0),\ell)} &= \bigcup_{i_1+ i_2+ i_3 \equiv \ell} (B_{i_1}\times B_{i_2}\times A_{i_3}), & 1 \leq \ell \leq 4^{n}, \\
\end{align*}
where the congruences are taken modulo $4^n$.

For each $\vec{b} \in \bodd_c$ and each $\ell \in [(2^{c-1})^{n}]$,
we define 
the sets $B_{(\vec{b},\ell)}$ similarly.
We show that the $(2^{c-1})^{n+1}$ sets $A_{(\vec{b},\ell)}$ form the desired partition of
$\beven_{c^{n+1}}$.
A similar proof shows that the $(2^{c-1})^{n+1}$ sets $B_{(\vec{b},\ell)}$ form the desired partition of $\bodd_{c^{n+1}}$.

Let  $\vec{x} \in \beven_{c^{n+1}}$.  
We partition the coordinates of 
$\vec{x}$ into $c$ bins of size $c^{n}$ and write $\vec{x} = (\vec{x}_1, \ldots, \vec{x}_c)$, where $\vec{x}_j$
is $\vec{x}$ restricted to bin~$j$.  
Let $b_j$ be the parity of $\vec{x}_j$ for each $j\in[c]$.
Since $\vec{x} \in \beven_{c^{n+1}}$, $(b_1, \ldots, b_c) \in \beven_c$. 
Since $\bigcup A_i = \beven_{c^{n}}$ and $\bigcup B_i = \bodd_{c^{n}}$, there exists an $\ell \in [(2^{c-1})^{n}]$
such that $\vec{x} \in A_{((b_1, \ldots, b_c), \ell)}$.
As the sets $A_{(\vec{b}, \ell)}$ are pairwise disjoint,
they  partition $\beven_{c^{n+1}}$.

Let $S$ be a subcube of $Q_{c^{n+1}}$ of dimension $c^{n+1} - c^{n} +1$.  Partition the coordinates of $S$ into $c$ bins
each of size $c^{n}$.  
Since $c^{n+1} - c^{n} +1 = (c-1)c^{n} + 1$, $\rS{j}$ has dimension at least~1 for every $j\in [c]$.
Since $c^{n+1} - c^{n} +1 = c(c^{n} - c^{n-1}) + 1$,  $\rS{j}$ has dimension at least
$c^{n} - c^{n-1} +1$ for some $j\in [c]$.  
W.l.o.g., $\rS{c}$ has dimension at least $c^{n} - c^{n-1} +1$.
Thus, $\rS{c}$ contains a vertex in every $A_i$ and $B_i$.
For  each $j \in [c-1]$, since $\rS{j}$ has dimension at least~1, 
$\rS{j}$ contains a vertex in each of $\beven_{c^n}$ and $\bodd_{c^n}$.
Thus, there exist two sequences, $k_1, k_2, \ldots, k_{c-1}$ and
$m_1, m_2, \ldots, m_{c-1}$ such that  $k_i, m_i \in [(2^{c-1})^{n}]$
for each $i \in [c-1]$ and 
$\rS{j}$ contains a vertex from $A_{k_j}$ and $B_{m_j}$ for each $j \in [c-1]$.

Let $(b_1, \ldots, b_c) \in \beven_c$ and let $\ell \in [(2^{c-1})^{n}]$. 
Let $E = \{j\in [c-1] : b_j = 0\}$ and let $F = \{j\in [c-1] : b_j = 1\}$.
For each $j\in [c-1]$, let
$D_j = A_{k_j}$ if $j \in E$ and $D_j = B_{m_j}$ if $j \in F$, where $A_{k_j}$ and $B_{m_j}$ 
are as defined above.
Let $i\in  [(2^{c-1})^{n}]$ satisfy 
$i + \sum_{j \in E} k_j + \sum_{j\in F} m_j  \equiv \ell \pmod{(2^{c-1})^{n}}.$
If $b_c = 0$, let $D_c = A_i$, otherwise, let $D_c = B_i$.
Since $\rS{j}$ contains a vertex from $D_j$ for each $j\in[c]$, 
$S$ contains a vertex from $(D_1\times \cdots \times D_c)$.
Thus, $S$ contains a vertex from $A_{(\vec{b}, \ell)}$ for
all $\vec{b}\in \beven_c$ and all $\ell \in [(2^{c-1})^{n}]$.
\hfill $\square$

\begin{theorem}\label{Q3d_d+1}
For each $d \geq 0$, there exist $4^{d+1}$ disjoint \BPSs\ for $\cQ(3^{d+1}, 3^d +1)$ 
with equal cardinalities
which partition the set of edges of $Q_{3^{d+1}}$.
\end{theorem}
The proof of Theorem~\ref{Q3d_d+1} (which we omit) is very similar to the proof of Theorem~\ref{Q4d_d+1},
except we use the following pairing strategies, which use rotating \binforms\ shifted by $s$:
\begin{align*}
\binPSP[n][3]{0}{s} & = (\bv \times \beven \times  \beven)^{R(s)} \cup   (\bodd \times  \bv \times \bodd)^{R(s)} \cup
(\beven \times  \bodd \times \bv)^{R(s)},  \\ \\
\binPSP[n][3]{1}{s} & = (\bv \times \beven \times  \bodd)^{R(s)} \cup   (\beven \times  \bv \times \beven)^{R(s)} \cup
(\bodd \times  \bodd \times \bv)^{R(s)},  \\ \\
\binPSP[n][3]{2}{s} & = (\bv \times \bodd \times  \beven)^{R(s)} \cup   (\beven \times  \bv \times \bodd)^{R(s)} \cup
(\bodd \times  \beven \times \bv)^{R(s)},  \\ \\
\binPSP[n][3]{3}{s} & = (\bv \times \bodd \times  \bodd)^{R(s)} \cup   (\bodd \times  \bv \times \beven)^{R(s)} \cup
(\beven \times  \beven \times \bv)^{R(s)}, 
\end{align*}
where, for example, 
\begin{align*}
(\beven \times  \bodd \times \bv)^{R(s)} = \{(\vec{x}, \vec{y}, \vec{z}) :\  & \vec{x} \in \beven, \vec{y}\in \bodd, 
\vec{z}\in\bvP{i}, \\
& \mbox{where $i = s + \Index(\vec{x})+\Index(\vec{y}) \pmod{m}$}\},
\end{align*}
where we assume that we have $m=4^d$ matchings $\bvP{i}$ (of equal cardinality) which partition
the edges of $Q_{3^d}$ and each $\bvP{i}$ is a \BPS\ for $\cQ(3^d, 3^{d-1}+1)$
in order to produce $4^{d+1}$ \BPSs\ for $\cQ(3^{d+1}, 3^d+1)$.

We use $\bvP{0}, \bvP{1}, \bvP{2}, \bvP{3}$ from the beginning of Section~\ref{Rotate}
for the case $d=0$.


\section{Conclusion}
\label{conclusion}

Let $\ps(n)$ be the smallest value of $k$ such that \B\ wins the positional game on
$\cQ(n,k)$ by {\em using a pairing strategy}.
We have proven the following upper bounds.
If $n\in\{4^{d+1} : d\in\naturals\}$, then $\ps(n) \leq \frac{n}{4} +1$.
If $n \in \{3^{d+1} : d \in \naturals\}\cup \{6 \cdot 4^d : d \in \naturals\} \cup \{9\cdot 4^ d : d\in\naturals\}$,
then $\ps(n) \leq \frac{n}{3} +1$.
In general, for all $n \geq 3$,
$\ps(n) \leq  \frac{3}{7}n +1$.
To obtain a lower bound on $\ps(n)$, we cite Proposition~9 in \cite{CPgames}, which implies that $\ps(n) > \ln(n)$.  Thus, there is a large gap between
the upper and lower bounds on $\ps(n)$ for most values of $n$.  
It would be interesting to improve any of these bounds.
With regards to small specific values of $n$,
because \M\
has a winning strategy for $\cQ(5,2)$ (see \cite{MBhypercube}) and $\cQ(2,1)$, we know 
that $\ps(3) = \ps(4) = 2$ and $\ps(5) = \ps(6) = 3$.  
It would be nice to also determine the exact values of, say, $\ps(7)$ and $\ps(8)$.

We note that there is no direct analogue to Theorems~\ref{Q4d_d+1} and \ref{Q3d_d+1}
for $\cQ(c^{d+1}, c^d +1)$ for $c \geq 5$ using our proof method.  Indeed,
Theorems~\ref{Q4d_d+1} and \ref{Q3d_d+1} rely on the \BPSs\
for $\cQ(4,2)$ and $\cQ(3,2)$ in order to create    
$\binPSP[n][4]{j}{s}$ and $\binPSP[n][3]{j}{s}$.  Since \M\ has     
a winning strategy for $\cQ(c,2)$ for all $c \geq 5$, there are no \BPSs\  
for $\cQ(c,2)$ from which we would create the \binforms\ for $\binPSP[n][c]{j}{s}$  for all $c \geq 5$.

As a final note, we mention that some of our results can be viewed as being related to a 
Tur\'an-type problem on $Q_n$.   Let $\ex(G, H)$ be the maximum number
of edges in a subgraph of $G$ which does not contain a copy of $H$.
In \cite{Erdos}, Erd\H{o}s discussed some problems that he believed
deserved more attention, including determining $\ex(Q_n, C_4)$,
which he conjectured to be $(\frac{1}{2} + o(1))|E(Q_n)|$.
Much work has been done related to determining $\ex(Q_n, C_{2t})$,
see for example, 
\cite{AlonEtAl}, 
\cite{AxenovichMartin}, 
\cite{BaloghEtAl}, 
\cite{Bialostocki}, 
\cite{BrassHarborthNienborg},
\cite{BrouwerEtAl},
\cite{Chung}, 
\cite{Conder}, 
\cite{Conlon}, 
\cite{FO:14-cycle, FurediOzkahya1}, 
\cite{ThomasonWagner}.

In \cite{AlonKrechSzabo}, Alon, Krech, and Szab\'o change the focus to 
studying $\ex(Q_n, Q_d)$. 
In particular, let $c(n, d)$ be 
the minimum number of edges that must
be deleted from $Q_n$ so that no copy of $Q_d$ remains,
and let $c_d = \lim_{n\to\infty} c(n,d)/ |E(Q_n)|$.
(For a study of $c(n,d)$ in a computer science context, see
\cite{GrahamHararyEtAl}.)
In their approach, Alon, Krech, and Szab\'o used a ``Ramsey-type framework," which involved
studying {\em $d$-polychromatic} colorings of the edges of $Q_n$, i.e.,
colorings in which every $d$-dimensional subcube of $Q_n$ contains
an edge from every color class.
They define $pc(n,d)$  to be the largest integer $p$ such that there exists
a $d$-polychromatic coloring of the edges of $Q_n$ in $p$ colors,
and $p_d =  \lim_{n\to\infty} pc(n,d)$.  They also define higher-dimensional
analogues, where the definition of $pc^{(\ell)}(n,d)$ is based on coloring each
$\ell$-dimensional subcube of $Q_n$ so that each $d$-dimensional
subcube contains an $\ell$-dimensional subcube of each color.  
Thus, $pc(n,d)$ is the special case $\ell = 1$.
They proved upper and lower bounds for $p_d$ for all $d \geq 1$ 
and that $p^{(0)}_d = d+1$ for all $d \geq 0$.
In~\cite{OffnerPoly}, Offner proved that $p_d$ equals the lower bound given by Alon, Krech, and Szab\'o.
Much work related to polychromatic colorings on the hypercube has been done,
for example, 
\cite{Chen},
\cite{GoldwasserEtAl},
\cite{HanOffner},
\cite{OffnerTType},
and 
\cite{OzkahyaStanton}.

We note that Theorems~\ref{Q3d_d+1} and \ref{Q4d_d+1} provide a $(3^d +1)$-polychromatic {\em proper} coloring of
$Q_{3^{d+1}}$ and a  $(4^d +1)$-polychromatic proper coloring of $Q_{4^{d+1}}$ for all $d \geq 0$, both using $4^{d+1}$ colors, i.e.,
each color class forms a matching.
It would be interesting to determine for which values of $n$ and $d$ there exists a  $d$-polychromatic {\em proper}
coloring of $Q_n$.

We also note that Lemma~\ref{PartitionEvenOdd} provides a $(c^n - c^{n-1} +1)$-polychromatic coloring of the {\em vertices}
of $Q_{c^n}$ using $(2^{c-1})^n$ colors and only vertices from $\beven_{c^n}$ (or $\bodd_{c^n}$).  
If we let $A_1, \ldots, A_{(2^{c-1})^n}$ be the partition of $\beven_{c^n}$ and
$B_1, \ldots, B_{(2^{c-1})^n}$ be the partition of $\bodd_{c^n}$, 
then $A_1 \cup B_1, \ldots, A_{(2^{c-1})^n}\cup B_{(2^{c-1})^n}$ works as a sort of 
$(c^n - c^{n-1} +1)$-polychromatic {\em double-coloring} of the vertices
of $Q_{c^n}$ using $(2^{c-1})^n$ colors, i.e., every $(c^n - c^{n-1} +1)$-dimensional subcube contains
{\em two} vertices from each color class.
It could be interesting to ask for which values of $n$, $d$, and $p$ do there exist $d$-polychromatic double-colorings of $Q_n$
using $p$ colors.

 \section{Acknowledgments}
The authors would like to thank Tam\'{a}s Lengyel 
for his valuable input and suggestions.
The authors would also like to thank the anonymous referee who
provided helpful comments and suggestions that improved the presentation of this paper.


\begin{thebibliography}{00}


\bibitem{AlonKrechSzabo} N. Alon, A. Krech, and T. Szab\'o,  Tur\'an's theorem in the hypercube,
{\it SIAM J.~Discrete Math.} {\bf 21} (2007), pp.~66-72.

\bibitem{AlonEtAl} N. Alon, R. Radoi\v{o}i\'c, B. Sudakov, and J. Vondra\'ak, A Ramsey-type result for the hypercube,
{\it J. Graph Theory} {\bf 53} (2006), 196--208.

\bibitem{AxenovichMartin} M. Axenovich and R. Martin, A note on short cycles in a hypercube,
{\it Discrete Math.} {\bf 306}, no.~18 (2006), pp.~2212--2218.


\bibitem{BaloghEtAl} J. Balogh, P. Hu, B. Lidick\'y, and H. Liu, Upper bounds on the size of 4- and 6-cycle-free
subgraphs of the hypercube, {\it European J. Combin.}, {\bf 35} (2014), 75--85.

\bibitem{TicTacToeBook} J.~Beck, {\it Combinatorial Games: Tic-Tac-Toe Theory}, 
Cambridge University Press, 2008.

\bibitem{BeckInevRandom} J.~Beck, {\it Inevitable Randomness in Discrete Mathematics},
 University Lecture Series, 49,  American Mathematical Society, Providence, RI, 
xii+250, 2009.


\bibitem{Bialostocki} A. Bialostocki, Some Ramsey type results regarding the graph of the $n$-cube,
{\it Ars Combinat.} 16-A (1983), pp.~39--48.

\bibitem{BrassHarborthNienborg} P. Brass, H. Harborth, and H. Nienborg, On the maximum number of edges in a 
$C_4$-free subgraph of $Q_n$, {\it J. Graph Theory} {\bf 19} (1995), 17--23.

\bibitem{BrouwerEtAl} A. E. Brouwer, I. J. Dejter, and C. Thomassen, Highly symmetric subgraphs of hypercubes,
{\it J. Algebraic Combin.} {\bf 2} (1993), 25--29.


\bibitem{Chen} E. Chen, Linear polychromatic colorings of hypercube faces,
{\it Electron. J. Combin.} {\bf 25} (1) (2018), \#P1.2.

\bibitem{Chung} F. Chung, Subgraphs of a hypercube containing no small even cycles, 
{\it J. Graph Theory} {\bf 16} (1992), 273--286.



\bibitem{Conder} M. Conder, Hexagon-free subgraphs of hypercubes, 
{\it J. Graph Theory} {\bf 17}, no.~4 (1993), pp.~477--479.


\bibitem{Conlon} D. Conlon, An extremal theorem in the hypercube, 
{\it Electron. J. Combin.} {\bf 17} (2010), \#R111




\bibitem{CPgames} A. Csernenszky, R. Martin, and A. Pluh\'ar, On the complexity of Chooser-Picker
positional games, {\it Integers} {\bf 11} (2011), G2, 16~pp.


\bibitem{Erdos} P. Erd\H{o}s, On some problems in graph theory, combinatorial analysis and combinatorial number theory,
{\it Graph Theory and Combinatorics}, B. Bollob\'as, ed., Academic Press (1984), 1--17.

\bibitem{ref:ES} P. Erd\H{o}s and J. L. Selfridge,  On a combinatorial game. {\it J. Combin.
Theory Ser. A}, {\bf 14} (1973), 298-301.

\bibitem{FO:14-cycle} Z. F\"uredi and L. \"Ozkahya, On 14-cycle-free subgraphs of the hypercube, 
{\it Combin. Probab. Comput.} {\bf 18} (2009), 725--729.

\bibitem{FurediOzkahya1} Z. F\"uredi and L. \"Ozkahya, On even-cycle-free subgraphs of the hypercube, 
{\it Electron. Notes  Discrete Math.} {\bf 34} (2009), 515--517.


\bibitem{GoldwasserEtAl} J. Goldwasser, B. Lidicky, R. Martin, D. Offner, J. Talbot, and M. Young,
Polychromatic colorings on the hypercube, {\it J. Comb.} {\bf 9}, no.~4 (2018), pp.~631--657.


\bibitem{GrahamHararyEtAl} N. Graham, F. Harary, M. Livingston, and Q. Stout, Subcube fault-tolerance in hypercubes,
{\it Inform. and Comput.} {\bf 102}, no.~2 (1993), 280--314.


\bibitem{HanOffner} E. Han, D. Offner, Linear $d$-polychromatic $Q_{d-1}$-colorings of the hypercube,
{\it Graphs Combin.} {\bf 34}, (2018) 791--801.

\bibitem{PosGamesBook} D. Hefetz, M. Krivelevich, M. Stojakovi\'c,  and
              T. Szab\'o,  {\it Positional Games},   Oberwolfach Seminars, 44, 
              Birkh\"auser/Springer, Basel, x+146, 2014.

\bibitem{MBhypercube} K. Kruczek and E. Sundberg, A Maker--Breaker game
on the boolean hypercube with subcubes as winning sets, {\it Integers} {\bf 18} (2018), G2, 21~pp.



\bibitem{OffnerPoly} D. Offner, Polychromatic colorings of subcubes of the hypercube, 
{\it SIAM J. Discrete Math.} {\bf 22}, no.~2 (2008), pp.~450--454.

\bibitem{OffnerTType} D. Offner, Some Tur\'an type results on the hypercube, 
{\it Discrete Math.} {\bf 309}, no.~9 (2009), pp.~2905--2912.

\bibitem{OzkahyaStanton} L. \"{O}zkahya and B. Stanton, On a covering problem in the hypercube,
{\it Graphs Combin.} {\bf 31}, no.~1 (2015), pp.~235--242.


\bibitem{ThomasonWagner} A. Thomason and P. Wagner, Bounding the size of square-free subgraphs of the hypercube,
{\it Discrete Math.} {\bf 309} (2009), 1730--1735.



\end{thebibliography}
\end{document}